\documentstyle[12pt]{article}
\newcommand{\const}{\mathop{\rm const}\limits}

\newcommand{\Law}{\mathop{\rm Law}\limits}

\newcommand{\Var}{\mathop{\rm Var}\limits}

\begin{document}

\begin{center}

{\bf CENTRAL LIMIT THEOREM AND} \par

\vspace{4mm}

{\bf EXPONENTIAL TAIL ESTIMATIONS IN}\par

\vspace{4mm}

{\bf MIXED (ANISOTROPIC) LEBESGUE SPACES }\par

\vspace{4mm}

 $ {\bf E.Ostrovsky^a, \ \ L.Sirota^b } $ \\

\vspace{4mm}

$ ^a $ Corresponding Author. Department of Mathematics and computer science, Bar-Ilan University, 84105, Ramat Gan, Israel.\\
\end{center}
E - mail: \ galo@list.ru \  eugostrovsky@list.ru\\
\begin{center}
$ ^b $  Department of Mathematics and computer science. Bar-Ilan University,
84105, Ramat Gan, Israel.\\

E - mail: \ sirota3@bezeqint.net\\

\vspace{3mm}
                    {\sc Abstract.}\\

 \end{center}

 \vspace{4mm}

 We  study the Central Limit Theorem (CLT) in the so-called  mixed (anisotropic)  Lebesgue-Riesz spaces
and tail behavior of normed sums of centered random independent variables (vectors)  with values in these spaces. \par

  \vspace{4mm}


\vspace{4mm}

 \par

\vspace{4mm}

\section{Notations. Statement of problem.}

\vspace{3mm}

{\bf 1.}  Let  $  (B, ||\cdot||B )  $  be separable Banach space and  $ \{  \xi_j   \}, \xi = \xi_1, \ j = 1,2,\ldots  $
be a sequence of centered in the weak sense: $ {\bf E} (\xi_i,b) = 0 \ \forall b \in B^* $  of independent identical distributed
(i.; i.d.)  random variables (r.v.) (or equally random vectors, with at the same abbreviation r.v.)  defined on some non-trivial
probability  space $  (\Omega = \{\omega\}, F, {\bf P})   $ with values  in the space  $ B. $  Denote

 $$
 S(n) = n^{-1/2} \sum_{j=1}^n \xi_j, \ n = 1,2,\ldots. \eqno(1.1)
 $$

 If we suppose that the r.v. $ \xi  $ has a weak second moment:

 $$
 \forall b \in B^* \ \Rightarrow  (Rb,b) := {\bf E } (\xi,b)^2 < \infty,  \eqno(1.2)
 $$
then the characteristical functional (more exactly,  the sequence of  characteristical functionals)

$$
\phi_{S(n)}(b) := {\bf E} e^{i \ (S(n), b) } \eqno(1.3)
$$
of $ S(n) $ converges as $ n \to \infty $ to the  characteristical functional of (weak, in general case)
Gaussian r.v. $ S = S(\infty)  $ with parameters $ (0,R):   $

$$
\lim_{n \to \infty} \phi_{S(n)}(b) = e^{ - 0,5 (Rb,b) }.
$$
  Symbolically:  $ S \sim N(0,R)  $ or $ \Law(S) = N(0,R).  $ The operator $ R = R_{S} $ is called
 the covariation operator, or variance of the r.v. $  S: $

$$
R = \Var(S);
$$
note that $  R = \Var(\xi). $\par

\vspace{4mm}

 We recall the classical definition of the CLT in the space $  B. $ \par

 \vspace{3mm}

 {\bf Definition 1.1.} {\it  We will say that the sequence $  \{ \xi_i  \} $  satisfies the CLT in the space $  B, $
 write:  $ \{\xi_j \} \in CLT = CLT(B) $ or   simple:  $ \xi \in CLT(B),  $  if the
 limit Gaussian r.v. $  S $  belongs to the  space $  B $ with probability one:  $  {\bf P} (S \in B) = 1 $ and the sequence
 of distributions $\Law(S(n))  $  converges weakly as $ n \to \infty $  to the distribution of the r.v. $ S = S(\infty):  $  }\par

$$
\lim_{n \to \infty} \Law(S(n)) = \Law(S). \eqno(1.4)
$$
 The  equality (1.4) imply that for any  continuous functional $ F: B \to R $

$$
\lim_{n \to \infty} {\bf  P} ( F(S(n)) < x  ) = {\bf P} (F(S) < x) \eqno(1.5)
$$
almost everywhere. \par

 In particular,

$$
\lim_{n \to \infty} {\bf  P} ( ||S(n)||B < x  ) = {\bf P} (||S||B < x), \ x > 0.
$$

 \vspace{3mm}

{\bf 2.} The problem of describing of necessary (sufficient) conditions  for the infinite - dimensional CLT in Banach space $  B $
has a long history; see, for instance, the  monographs \cite{Araujo1} - \cite{Ostrovsky1} and articles
\cite{Garling1} - \cite{Zinn1} ; see also reference therein.\par
 The applications of considered theorem in statistics and method Monte-Carlo  see, e.g. in
\cite{Frolov1} -  \cite{Ostrovsky303}.\par

\vspace{3mm}

{\bf 3.} The cornerstone of  this problem is to establish the {\it weak compactness } of the distributions generated
in the space $  B  $ by the sequence $  \{  S(n) \}: $

$$
\nu_n(D) = {\bf P} ( S(n) \in D),
$$
where $  D  $ is Borelian set in $  B; $ see \cite{Prokhorov1}; \cite{Billingsley1},  \cite{Billingsley2}. \par

\vspace{3mm}

{\bf 4.} We recall here the definition of
 the so-called anisotropic Lebesgue (Lebesgue-Riesz) spaces, which appeared in the famous article of
 Benedek A. and Panzone  R. \cite{Benedek1}.  More detail information about this
spaces with described applications see in the books  of  Besov O.V., Il’in V.P., Nikol’skii S.M.
\cite{Besov1}, chapter 1,2; Leoni G. \cite{Leoni1}, chapter 11; \cite{Lieb1}, chapter 6. \par

\vspace{3mm}

  Let $ (X_k,A_k,\mu_k), \ k = 1,2,\ldots,l $ be measurable spaces with sigma-finite {\it separable}
non - trivial measures $ \mu_k. $ The separability  denotes that  the metric space
$ A_k  $  relative the distance

$$
\rho_k(D_1, D_2) = \mu_k(D_1 \Delta D_2) = \mu_k(D_1 \setminus D_2) + \mu_k(D_2 \setminus D_1) \eqno(1.6)
$$
is separable.\par
Let also $ p = (p_1, p_2, . . . , p_l) $ be $ l- $ dimensional vector such that
$ 1 \le p_j < \infty.$ \par

 Recall that the anisotropic  (mixed) Lebesgue - Riesz space $ L_{ \vec{p}} $ consists on all the  total measurable
real valued function  $ f = f(x_1,x_2,\ldots, x_l) = f( \vec{x} ) $

$$
f:  \otimes_{k=1}^l X_k \to R
$$

with finite norm $ |f|_{ \vec{p} } \stackrel{def}{=} $

$$
\left( \int_{X_l} \mu_l(dx_l) \left( \int_{X_{l-1}} \mu_{l-1}(dx_{l-1}) \ldots \left( \int_{X_1}
 |f(\vec{x})|^{p_1} \mu(dx_1) \right)^{p_2/p_1 }  \ \right)^{p_3/p_2} \ldots   \right)^{1/p_l}. \eqno(1.7)
$$

 In particular, for the r.v. $ \xi  $

$$
  |\xi|_p =  \left[ {\bf E} |\xi|^p \right]^{1/p}, \ p \ge 1.
$$

 Note that in general case $ |f|_{p_1,p_2} \ne |f|_{p_2,p_1}, $
but $ |f|_{p,p} = |f|_p. $ \par

 Observe also that if $ f(x_1, x_2) = g_1(x_1) \cdot g_2(x_2) $ (condition of factorization), then
$ |f|_{p_1,p_2} = |g_1|_{p_1} \cdot |g_2|_{p_2}, $ (formula of factorization). \par

 Note that under conditions separability (1.6) of measures $ \{  \mu_k \} $   this spaces are also  separable and Banach spaces. \par

  These spaces arises in the Theory of Approximation, Functional Analysis, theory of Partial Differential Equations,
theory of Random Processes etc. \par

\vspace{3mm}

{\bf 5.} Let for example $  l = 2; $ we agree to rewrite for clarity the expression for $ |f|_{p_1, p_2}  $ as follows:

$$
|f|_{p_1, p_2} := |  f|_{p_1, X_1; p_2, X_2}.
$$
 Analogously,

$$
|f|_{p_1, p_2,p_3} = |  f|_{p_1, X_1; p_2, X_2; p_3, X_3}. \eqno(1.8)
$$

  Let us give an example. Let $ \eta = \eta(x, \omega )  $ be bi - measurable  random field, $ (X = \{x\}, A, \mu) $ be measurable space,
  $  p = \const \in [1,\infty). $ As long as the expectation $ {\bf E } $ is also an integral, we deduce

  $$
{\bf E} |\eta(\cdot, \cdot)|^p_{p,X} = {\bf E} \int_X |\eta(x, \cdot)|^p \ \mu(dx)    =
  $$

$$
 \int_X  {\bf E}  |\eta(x, \cdot)|^p \ \mu(dx)  = \int_X \mu(dx) \left[ \int_{\Omega} |\eta(x, \omega)|^p \ {\bf P}(d \omega)  \right];
$$

$$
| \ |\eta|_{p,X} \ |_{p, \Omega} = \left\{ \int_X \mu(dx) \left[ \int_{\Omega} |\eta(x, \omega)|^p \ {\bf P}(d \omega)  \right] \right\}^{1/p} =
|  \ \eta(\cdot, \cdot)  \  |_{1,X; p,\Omega}. \eqno(1.9)
$$

\vspace{3mm}

{\bf 6.}  Constants of Rosenthal - Dharmadhikari - Jogdeo - ...\par

 Let  $  p = \const \ge 1, \hspace{4mm}  \{ \zeta_k \} $ be a sequence of numerical centered, i.; i.d. r.v.  with finite $ p^{th} $ moment
 $ | \zeta|_p < \infty. $  The following constant,  more precisely, function on $ p, $ is called
 constants of Rosenthal - Dharmadhikari - Jogdeo - Johnson - Schechtman - Zinn - Latala - Ibragimov - Pinelis - Sharachmedov - Talagrand - Utev...:

$$
K_R(p) \stackrel{=} \sup_{n \ge 1} \sup_{ \{\zeta_k\} } \left[ \frac{|n^{-1/2} \sum_{k=1}^n \zeta_k|_p}{|\zeta_1|_p} \right].  \eqno(1.10)
$$
 We will  use  the following ultimate up to an error value $ 0.5\cdot 10^{-5} $  estimate  for $ K_R(p), $ see \cite{Ostrovsky502} and reference therein:

 $$
 K_R(p) \le \frac{C_R \ p}{ e \cdot \log p}, \hspace{5mm}  C_R = \const := 1.77638.  \eqno(1.11)
 $$
 Note that for the symmetrical distributed r.v. $ \zeta_k $ the constant $  C_R $ may be reduced  up to a value $ 1.53572.$\par

\vspace{3mm}

{\bf 7. } The Law of Large Numbers  (LLN) on Banach spaces.\\

\vspace{2mm}

It is known, see, e.g. \cite{Fortet1},  that the sequence $  \{ \xi_k \} $ of mean zero, i.,i.d. random vectors   satisfies
{\it Strong} LLN in the separable Banach space $ (B, ||\cdot||) $ if and only if $ {\bf E} ||\xi_1|| < \infty. $ \\

\vspace{3mm}

\section{ Ordinary Lebesgue - Riesz spaces.  }

\vspace{3mm}
  We study in this section the CLT and exponential tail estimates  for the normed sums of centered, i., i.d. in the Banach space
 $ L_p, \ 2 \le p < \infty. $ Notice that the case $ 1 \le p < 2 $ is considered in  \cite{Gine2}; see also reference therein.\par

  In detail,  $ \xi(x) =  \xi(x,\omega) $
 be bi - measurable centered  random field, $ (X = \{x\}, A, \mu) $ be  measurable space with separable sigma - finite measure $ \mu, $
 and $ \{ \xi_k(x,\omega)  \} = \{ \xi_k(x)   \} $ be independent copies of $ \xi = \xi(x) := \xi(x, \omega), $

$$
S_n(x) = n^{-1/2} \sum_{k=1}^n \xi_k(x). \eqno(2.1)
$$

\vspace{4mm}

{\bf Theorem 2.1.} {\it  If   }
$$
 {\bf E} |\xi|_{p,X}^p = {\bf E} \int_X |\xi(x)|^p \ \mu(dx) = \int_X {\bf E} |\xi(x)|^p \mu(dx) < \infty,
$$
{\it then  } $  \{\xi_k(\cdot) \} \in CLT  $ {\it in the space  } $  L_p(X). $\par

\vspace{4mm}

{\bf Theorem 2.2.} {\it  Let $ m = \const \ge 1 $ be number, not necessary to be integer,
 for which }

$$
 {\bf E} | \xi(\cdot)|_{pm,X}^{pm} < \infty. \eqno(2.2.)
$$

{\it or equally }

$$
 |\xi(\cdot) \ |_{p,X; mp, \Omega}  < \infty.  \eqno(2.2a)
$$

 {\it or equally }

$$
{\bf E }  \int_X |\xi(x)|^{pm} \ \mu(dx) = \int_X {\bf E } |\xi(x)|^{pm} \ \mu(dx) < \infty.  \eqno(2.2b)
$$

{\it Then }

$$
\sup_n {\bf E} | S_n(\cdot)|_{p,X}^{pm} \le K_R^{pm} (pm) \  \left[ \int_X \left[ {\bf E} | \xi(x)|^{pm}  \right]^{1/m} \ \mu(dx) \right]^m, \eqno(2.3)
$$

 {\it or equally}

 $$
 \sup_n | S_n|_{p, X; mp,\Omega} \le K_R(pm) \ |\xi|_{pm, \Omega; p,X}. \eqno(2.3a).
 $$

\vspace{3mm}

{\bf Proofs.}\\

\vspace{3mm}

{\bf 0.} The proposition of theorem 2.1. is known; see, e.g. \cite{Ledoux1}, chapter 10, p. 273 - 278;
 \cite{Grenander1}, chapter 8; but we  intend to
prove it through the inequality (2.3) and we add the non-asymptotical tail estimates for normed sums of independent
mean zero random fields.  \par

 Moreover, as long as the space $  L_p(X) $ obeys a type 2, the {\it necessary and sufficient } condition for CLT in the space
$  L_p(X) $  has a view:

$$
{\bf E} \xi = 0, \hspace{5mm} {\bf E} ||\xi||^2 < \infty,
$$
see also \cite{Ledoux1}, chapter 10, p.  281. \par

  Note that in the arbitrary separable Banach space $ B $ the condition

$$
\lim_{t \to \infty} t^2 {\bf P} ( ||\xi|| > t) = 0
$$
is {\it necessary} for CLT in the space $  B, $ see \cite{Ledoux1}, chapter 10, p.  274. \par
In particular, $ \forall q \in (0,2)  \Rightarrow {\bf E} ||\xi||^q < \infty. $\par

\vspace{3mm}

{\bf Remark 2.1.} The mixed Lebesgue space $ L_{\vec{p}} $ does not have in general case the type 2. It suffices to
consider the space $ L_{3, 3/2, 5, 3/2}[0,1]^4.  $ \par

\vspace{3mm}

{\bf 1. Auxiliary fact. } \\

\vspace{3mm}

{\bf Lemma 2.1.} {\it Under condition 2.2  }

$$
{\bf E} \left[ \int_X |\xi(x)|^p \ \mu(dx)  \right]^m \le  \left\{ \int_X \left[ {\bf E} |\xi(x)|^{pm} \right]^{1/m} \  \mu(dx) \right\}^m, \eqno(2.4)
$$
{\it  or equally}

$$
|\xi|_{p,X; pm, \Omega} \le |\xi|_{pm, \Omega; p,X}. \eqno(2.4a)
$$

\vspace{3mm}

{\bf Proof of Lemma 2.1.}

\vspace{3mm}

 Denote $  \eta(x) = |\xi(x)|^p.  $ Note that the  space $ L_m(\Omega)  $ is a Banach space; and we can apply the
generalized triangle (Minkowsky) inequality:

$$
\left[ {\bf E}  \left(  \int_X |\xi(x)|^p \ \mu(dx)  \right)^m \right]^{1/m}  =
$$

$$
\left[ {\bf E}  \left(  \int_X \eta(x) \ \mu(dx)  \right)^m \right]^{1/m}  = \left| \ \int_X \eta(x) \ \mu(dx)   \  \right|_{m, \Omega} \le
$$

$$
 \ \int_X  |\eta(x)|_{m, \Omega} \ \mu(dx) = \int_X \sqrt[m] {{\bf E} \eta^m(x) } \ \mu(dx)  =
\int_X \sqrt[m] {{\bf E} |\xi|^{pm}  (x) } \ \mu(dx), \eqno(2.5)
$$
which is equivalent to the assertion (2.4). \par
 We  used theorem of Fubini.

\vspace{4mm}

{\bf 2. Proof of theorem 2.2.}  It is sufficient to apply the assertion of Lemma 2.1. to the random field
$  S_n(x)  $ instead $ \xi(x): $

$$
{\bf E} \left[ \int_X |S_n(x)|^p \ \mu(dx)  \right]^m \le
\left\{ \int_X  \left[{\bf E}|S_n(x)|^{pm} \ \right]^{1/m} \mu(dx) \right\}^m.  \eqno(2.6)
$$

 We obtain by means of Rosenthal's inequality (recall that we  consider here the case only when $  p \ge 2:) $

$$
{\bf E} |S_n(x)|^{pm} \le K_R^{p m} (p m) \ {\bf E} |\xi(x)|^{pm} =  K_R^{p m}(p m) \cdot |  \xi(x)|^{pm}_{pm, \Omega}.
$$
 It remains to substitute into (2.6). \par

\vspace{3mm}

{\bf 3. Proof of theorem 2.1.}
 Let $  {\bf E} |\xi|_{p,X}^p < \infty. $  As long as the Banach space $ L_p(X) $  is separable
and the function $  y \to |y|^p  $ satisfies the $ \Delta_2  $ condition, there exists a linear compact operator
$  U: \ L_p(X) \to L_p(X)  $ such that

$$
{\bf P} \left(U^{-1} \xi \in L_p(X) \right) = 1 \eqno(2.7)
$$
and moreover

$$
{\bf E} | U^{-1} \xi |^p_{p,X} < \infty. \eqno(2.8)
$$
 \cite{Ostrovsky2}; see also \cite{Buldygin1}, \cite{Ostrovsky603}.\par

  Let us consider the sequence of r.v. in the space $ L_p(X): \ \eta_k(x) = U^{-1} [\xi_k](x); $ it is also a
 sequence of i., i.d. r.v. in the space $  L_p(X), $  and we can apply the inequality (2.3) tacking the value $  m = 1: $

$$
\sup_n {\bf E} | U^{-1}[ S_n]|_{p,X}^{p} \le K_R^{p} (p) \ {\bf E} | U^{-1}[\xi]|_{p,X}^{p} = C(p) < \infty. \eqno(2.9)
$$
 We get using Tchebychev's  inequality

 $$
\sup_n {\bf P} \left( |U^{-1}[ S_n]|_{p,X} > Z   \right) \le C(p)/Z^p < \epsilon,\eqno(2.10)
 $$
for sufficiently greatest values $ Z = Z(\epsilon), \ \epsilon \in (0,1). $ \par
  Denote by $  W  = W(Z) $ the set

$$
W = \{ f: f \in L_p(X),  |U^{-1}[f]|_{p,X}  \le Z \}. \eqno(2.11)
$$
 Since the operator $ U $ is compact, the set $  W  = W(Z) $ is compact set in the space $  L_p(X). $ It follows from inequality
(2.8) that

$$
\sup_n {\bf P} \left( S(n) \notin W(Z) \right) \le \epsilon.
$$
  Thus, the sequence $ \{  S_n \}  $ satisfies the famous Prokhorov's criterion \cite{Prokhorov1} for  weak compactness of
the family of distributions in the separable metric spaces. \par
 This completes the proof of theorem 2.1. \par

\vspace{4mm}
{\bf Examples.} \\

 Suppose

 $$
 {\bf P}  \left( |\xi|^p_{p,X} > x \right) \le e^{-x^Q }, \ x > 0, \ Q = \const > 0;
 $$
then  of course

$$
\forall m > 0 \ \Rightarrow {\bf E} |\xi|_{p,X}^{ p m} < \infty;
$$
therefore $ \{ \xi_k \}  $  satisfies CLT in the space $  L_p(X) $  and

$$
\sup_n {\bf P} \left( |S_n|^p_{p,X} > x  \right) \le \exp \left( - C(Q,p) \ x^{Q/(Q+1)   }  \right),  C(Q,p) > 0.
$$

 More generally, if

$$
 {\bf P}  \left( |\xi|^p_{p,X} > x \right) \le e^{-x^Q_1 \ (\log x)^{-Q_2} }, \ Q_1 = \const > 0, Q_2 = \const, \ x \ge e,\eqno(2.12)
 $$
then for some positive value $ C_3 = C_3(Q_1, Q_2,p) $ and $  x > e $

$$
\sup_n {\bf P} \left( |S_n|^p_{p,X} > x  \right) \le
$$

$$
\exp \left( - C_3(Q_1, Q_2,p) \ x^{Q_1/(Q_1+1)} \ ( \log x )^{( - Q_2 -  Q_1(Q_1-1)  )/( Q_1 + 1 ) }  \right).\eqno(2.13)
$$

 We used some estimations from the monograph  \cite{Ostrovsky1}, chapter 2, section 3, p. 55 - 57. \par

\vspace{3mm}
\section{ Mixed Lebesgue - Riesz spaces. Main results.}
\vspace{3mm}

   Let us return to the CLT in the space  $ L_{ \vec{p}}, \ \vec{p} = \{ p_k  \}, k=1,,2,\ldots,l, \ l  \ge 2; $
 where $  1 \le p_k < \infty, $ described below. \par

 Define

$$
\overline{p} : = \max(p_1,p_2, \ldots, p_l).
$$
{\it and suppose  everywhere further}

$$
 \overline{p} \ge 2. \eqno(3.0)
$$

\vspace{4mm}

{\bf Theorem 3.1.} {\it  If in addition }

$$
 \left[ {\bf E} |\xi(\vec{x})|^{\overline{p}} \right]^{1/\overline{p}} \in L_{ \vec{p} }, \eqno(3.1)
$$
{\it then  } $  \{\xi_k(\cdot) \} \in CLT  $ {\it in the space  } $  L_{\vec{p}}. $\par

\vspace{4mm}

{\bf Theorem 3.2.} {\it  Let $ m = \const \ge 1  $ be (not necessary to be integer) number for which }

$$
 \left[ {\bf E} |\xi(\vec{x})|^{m \ \overline{p}} \right]^{1/(m \ \overline{p}) } \in L_{ \vec{p} }. \eqno(3.2)
$$

{\it Then }

$$
\sup_n  \left[ {\bf E}  | S_n(\cdot)|_{\vec{p}}^{ m \ \overline{p}} \right]^{ 1/(m \ \overline{p}) }  \le K_R(\overline{p} \ m)  \times
  \left| \left[ {\bf E} |\xi(\vec{x})|^{m \ \overline{p}} \right]^{1/(m \overline{p})} \right|_{\vec{p}}.\eqno(3.3)
$$

\vspace{3mm}

{\bf Proofs.} \\
\vspace{3mm}
{\bf 1. Auxiliary fact: permutation inequality.} \par

\vspace{3mm}

 We will use the so-called {\it permutation inequality} in the terminology of an article \cite{Adams1};  see also \cite{Besov1}, chapter 1,  p. 24 - 26.
Indeed, let $ (Z, B, \nu) $  be another measurable space  and $ \phi:  (\vec{X},Z) = \vec{X} \otimes Z    \to R $ be measurable function. In what follows
$  \vec{X} = \otimes_k X_k.   $ Let also $  r = \const \ge \overline{p}. $ It is true the following inequality (in our notations):

$$
|\phi|_{\vec{p}, \vec{X}; r, Z } \le |\phi|_{r, Z;  \vec{p}, \vec{X}}. \eqno(3.4)
$$
 In what follows $  Z = \Omega, \ \nu = {\bf P}.  $

\vspace{3mm}

{\bf 2. Auxiliary inequality.} \par

\vspace{3mm}

 It follows from permutation inequality (3.4) that

 $$
 \sqrt[m \overline{p}] {{\bf E} |\xi|_{\vec{p}}^{m \overline{p}}}  \le
 \left| \sqrt[m \overline{p}] { {\bf E} |\xi|^{m \overline{p}}} \right|_{\vec{p}}, \ m = \const \ge 1. \eqno(3.5)
 $$

\vspace{3mm}
{\bf 3. Proof of theorem 3.2.}\\
\vspace{3mm}

 We deduce applying  the inequality (3.5) for the random field $ S_n  $ and using the Rosenthal's inequality:

$$
 \sqrt[m \overline{p}] {{\bf E} |S_n|_{\vec{p}}^{m \overline{p}}}  \le K_R( m \ \overline{p} ) \cdot
 \left| \sqrt[m \overline{p}] { {\bf E} |\xi|^{m \overline{p}}} \right|_{\vec{p}}, \ m = \const \ge 1, \eqno(3.6)
 $$
which is equivalent to the assertion of theorem 3.2.\par

\vspace{3mm}

{\bf Proof of theorem 3.1} is alike ones in theorem 2.1. Let

$$
 \left| \left[ {\bf E} |\xi(\vec{x})|^{\overline{p}} \right]^{1/\overline{p}} \right|_{ \vec{p} } < \infty.  \eqno(3.7)
$$

 We use the proposition of theorem 3.2 with the value $  m = 1: $

$$
\sup_n {\bf E}  \left[ | S_n(\cdot)|_{\vec{p}}^{ \overline{p}} \right]^{ 1/ \overline{p} }  \le K_R(\overline{p} )  \times
  \left| \left[ {\bf E} |\xi(\vec{x})|^{\overline{p}} \right]^{1/\overline{p}} \right| < \infty.  \eqno(3.8)
$$

 As long as the Banach space $ L_{\vec{p}} = L_{\vec{p}, \vec{X}} $  is separable
and the function $  z \to |z|^{\overline{p}}, \ z \in R  $ satisfies the $ \Delta_2  $ condition, there exists a linear compact (non-random!)
 operator $  U: \ L_{\vec{p}} \to L_{\vec{p}} $ such that

$$
{\bf P} \left(U^{-1} \xi \in L_{\vec{p}} \right) = 1
$$
and moreover

$$
{\bf E} | U^{-1} \xi |^{\overline{p}}_{\vec{p}}  < \infty.
$$
 \cite{Ostrovsky2}; see also \cite{Buldygin1}, \cite{Ostrovsky603}.\par

  Let us consider the sequence of r.v. in the space $ L_{\vec{p}}: \ \eta_k(x) = U^{-1} [\xi_k](x); $ it is also a
 sequence of i., i.d. r.v. in the space $ L_{\vec{p}}, $  and we can apply the inequality of theorem 3.2:

$$
\sup_n {\bf E} | U^{-1}[ S_n]|_{\vec{p}}^{\overline{p}} \le K_R^{\overline{p}} (\overline{p}) \
{\bf E} | U^{-1}[\xi]|_{\vec{p}}^{\overline{p}} = C(\vec{p}) < \infty.  \eqno(3.9)
$$
 We get using Tchebychev's  inequality

 $$
\sup_n {\bf P} \left( |U^{-1}[ S_n]|_{\vec{p}} > Z   \right) \le C(\vec{p})/Z^{\overline{p}} < \epsilon,
 $$
for sufficiently greatest values $ Z = Z(\epsilon), \ \epsilon \in (0,1). $ \par
  Denote by $  W  = W(Z) $ the set

$$
W = \{ f: f \in L_{\vec{p}}, \  |U^{-1}[f]|_{\vec{p}}  \le Z \}.
$$
 Since the operator $ U $ is compact, the set $  W  = W(Z) $ is compact set in the space $  L_{\vec{p}, \vec{X}}. $ It follows from inequality
of lemma 3.1 that

$$
\sup_n {\bf P} \left( S(n) \notin W(Z) \right) \le \epsilon.  \eqno(3.10)
$$
  Thus, the sequence $ \{  S_n \}  $ satisfies the famous Prokhorov's criterion \cite{Prokhorov1} for  weak compactness of
the family of distributions in the separable metric spaces. \par
 This completes the proof of theorem 3.1. \par

\vspace{5mm}

\section{ Concluding remarks.}
\vspace{3mm}
{\bf A. Necessity of our conditions.}
\vspace{3mm}
 In the case when $ p_1 = p_2 = \ldots = p_l = 2  $ the space $ L_{\vec{p}} =   L_{\vec{p}, \vec{X}} $
is ordinary Hilbert space, for which the necessary and sufficient condition for CLT coincides with our
condition theorem 3.1:

$$
{\bf E} ||\xi||^2 < \infty.
$$

\vspace{3mm}
{\bf B. CLT in mixed Sobolev's spaces.}\\
\vspace{3mm}

 The method presented here may be used by investigation of Central Limit Theorem in the so-called
{\it mixed Sobolev's spaces} $ W^A_{\vec{p}}. $ In detail, let $ (Y_k, B_k, \zeta_k), \ l = 1,2,\ldots,l  $ be again
measurable spaces with separable sigma - finite measures $ \zeta_k. $
 Let $  A  $ be closed unbounded operator acting from the space $ W_{\vec{p}, \vec{X}} $ into
 the space $ W_{\vec{p}, \vec{Y}}, $  for instance, differential operator,
 Laplace's operator or its power, may be fractional, for instance:

 $$
 A[u](x,y) = \frac{D^{(\vec{q})}u(x) - D^{(\vec{q})}u(y)}{|x-y|^{\beta}}, \ x,y \in R^d, \ x \ne y, \
 $$

 $$
 \zeta(G) = \int \int_{G} \frac{dx dy}{|x-y|^{\alpha}}, \ \alpha,\beta = \const \in [0,1],  \alpha + \beta p < d, \ G \subset R^{2d},
 $$

$$
\vec{q} = \{q_1, q_2,  \ldots, q_d \}, \ q_s = 0,1, \ldots,
$$

$$
D^{(\vec{q})}u(x) = \frac{\partial^{q_1}}{\partial x_1^{q_1}} \frac{\partial^{q_2}}{\partial x_2^{q_2}} \ldots    \frac{\partial^{q_d}}{\partial x_d^{q_d}}u(x).
$$
  The norm in this space may be defined as follows (up to closure):

 $$
|f|W^A_{\vec{p}}\stackrel{def}{=} \max \left( |f|_{\vec{p}}, |Af|_{\vec{p}} \right).
 $$

 Analogously to the proof of theorems 3.1, 3.2 may be obtained the following results.\par

\vspace{3mm}

{\bf Theorem 4.1.} {\it If as before } $ \overline{p} \ge 2.$
{\it  and  }

$$
 \left[ {\bf E} |\xi(\vec{x})|^{\overline{p}} \right]^{1/\overline{p}} \in L_{ \vec{p} }, \hspace{5mm}
 \left[ {\bf E} | A[\xi](\vec{x})|^{\overline{p}} \right]^{1/\overline{p}} \in L_{ \vec{p} }, \eqno(4.1)
$$
{\it then  } $  \{\xi_k(\cdot) \} \in CLT  $ {\it in the space  } $  W^A_{\vec{p}}. $\par

\vspace{4mm}

{\bf Theorem 4.2.} {\it  Let $ m = \const \ge 1  $ be (not necessary to be integer) number for which }

$$
 \left[ {\bf E} |\xi(\vec{x})|^{m \ \overline{p}} \right]^{1/(m \ \overline{p}) } \in L_{ \vec{p} }, \hspace{5mm}
 \left[ {\bf E} |A[\xi](\vec{x})|^{m \ \overline{p}} \right]^{1/(m \ \overline{p}) } \in L_{ \vec{p} }. \eqno(4.2)
$$

{\it Then }

$$
\sup_n \left[ {\bf E} |A[ S_n](\cdot)|_{\vec{p}}^{ m \ \overline{p}} \right]^{ 1/(m \ \overline{p}) }  \le K_R(\overline{p} \ m)  \times
  \left| \left[ {\bf E} |A[\xi](\vec{x})|^{m \ \overline{p}} \right]^{1/(m \overline{p})} \right|_{\vec{p}}.\eqno(4.3)
$$

\vspace{3mm}
{\bf B. CLT for dependent r.v. in anisotropic spaces.}\\
\vspace{3mm}

{\it We refuse in this section on the assumption about independence of random vectors} $ \{ \xi_k(\cdot)  \}.  $ \\

\vspace{3mm}

{\bf  Martingale case.} \\

\vspace{3mm}

 We suppose as before that  $ \{ \xi_k(\cdot)  \}  $ are mean zero and form a strictly stationary sequence,
 $  \overline{p} \ge 2. $
 Assume in addition that $ \{ \xi_k(\cdot)  \}  $ form a martingale difference sequence
 relative certain filtration $ \{  F(k) \},  \ F(0) = \{ \emptyset, \Omega\},  $
 $$
  {\bf E} \xi_k/F(k) = \xi_k, \hspace{5mm}  {\bf E} \xi_k/F(k-1) = 0, \ k= 1,2,\ldots.
 $$

 Then the proposition  of theorem 3.1 remains true; the estimate (3.3) of theorem 3.2 is also true  up to multiplicative
absolute constant. \par
 Actually, the convergence of correspondent characteristical functionals follows from the ordinary CLT for martingales,
 see in the classical monograph of  Hall P., Heyde C.C.
  \cite{Hall1}, chapter 2; the Rosenthal's constant for the sums  of martingale differences with at the same up to
  multiplicative constant coefficient is obtained by A.Osekowski \cite{Osekowski1}, \cite{Osekowski2}.
 See also \cite{Ostrovsky3}.\par

\vspace{3mm}

{\bf  Mixingale case.} \\

\vspace{3mm}

 We suppose again  that  $ \{ \xi_k(\cdot)  \}  $ are mean zero and form a strictly stationary sequence,
 $  \overline{p} \ge 2. $  This  sequence is said to be {\it mixingale, } in the terminology of the book
\cite{Hall1}, if it satisfies this or that mixing condition.\par

 We  consider here only the  superstrong mixingale. Recall that the superstrong, or $  \beta = \beta(F_1, F_2) $
index between two sigma - algebras  is defined as follows:

$$
\beta(F_1, F_2) = \sup_{A \in F_1, B \in F_2, {\bf P}(A) {\bf P}(B) > 0 } \left| \frac{{\bf P}(AB) - {\bf P}(A) {\bf P}(B)}{{\bf P}(A) {\bf P}(B)} \right|.
$$

 Denote

 $$
 F_{-\infty}^0  = \sigma(\xi_s, \ s \le 0),  \hspace{5mm} F_n^{\infty} = \sigma(\xi_s, \ s \ge n), \eqno(4.4)
 $$

$$
\beta(n) = \beta \left(F_{-\infty}^0 , F_n^{\infty} \right),
$$

 The  sequence $ \{\xi_k \} $ is said to be {\it superstrong mixingale, } if $ \lim_{n \to \infty} \beta(n) = 0. $ \par
This notion  with some applications was introduced and investigated by B.S.Nachapetyan and R.Filips \cite{Nachapetyan1}.
See also   \cite{Ostrovsky3}, \cite{Ostrovsky1}, p. 84 - 90. \par

\vspace{3mm}

 Introduce the so-called mixingale Rosenthal coefficient:

 $$
 K_M(m) =  m \ \left[  \sum_{k=1}^{\infty} \beta(k) \ (k+1)^{ (m -2)/2  }   \right]^{1/m}, \ m \ge 1.\eqno(4.5)
 $$

B.S.Nachapetyan in \cite{Nachapetyan1} proved that for the superstrong centered  strong stationary strong mixingale
sequence $  \{ \eta_k \} $ with $ K_M(m) < \infty $

$$
\sup_{n \ge 1}  \left| n^{-1/2} \sum_{k=1}^n \eta_k \right|_m \le C \cdot K_M(m) \cdot |\eta_1|_m, \eqno(4.6)
$$
so that the "constant" $ K_M(m) $ play at the same role for mixingale as the Rosenthal  constant  for independent variables.\par

 As a consequence:  theorems 3.1 and 3.2  remains true for strong mixingale sequence $ \{  \xi_k \}: $
theorem 3.1 under conditions: $ K_M(1) < \infty $  for theorem 3.1 and $ K_M(m \ \overline{p}) < \infty $  for the theorem 3.2
with replacing $  K_R( m \ \overline{p} ) $ on the expression $  K_M( m \ \overline{p} ). $  \par


\vspace{6mm}
\begin{thebibliography}{99}

\vspace{4mm}

\bibitem{Adams1}
{\sc Adams R.A.} {\it  Anisotropic Sobolev Inequalities.}
Casopic pro Pestovani Matematiky, (Prague),  No. 3, 267—279.

\bibitem{Benedek1}
{\sc Benedek A. and Panzone  R.} {\it The space $  L_p $ with mixed norm.} Duke Math. J., {\bf 28}, (1961),  301 - 324.

\bibitem{Besov1}
{\sc Besov O.V., Il’in V.P., Nikol’skii S.M. } {\it Integral representation of functions
and imbedding theorems.} Vol.1; Scripta Series in Math., V.H.Winston
and Sons, (1979), New York, Toronto, Ontario, London.

\bibitem{Buldygin1}
{\sc Buldygin V.V.} (1984). {\it Supports of probabilistic measures in separable Banach
spaces.} Theory Probab. Appl. 29 v.3, pp. 528 - 532, (in Russian).


\bibitem{Hall1}
{\sc Hall P., Heyde C.C. } {\it Martingale Limit Theory and Applications.}  Academic
Press, New York. (1980)

\bibitem{Nachapetyan1}
{\sc Nachapetyan B.S.} {\it On the certain criterion of weak dependence. } Probab. Theory Appl., (1980),
{\bf 2, } V. 26, 374 - 381.

\bibitem{Leoni1}
{\sc Leoni G. } {\it A first Course in Sobolev Spaces.} Graduate Studies in Mathematics,
v. 105, AMS, Providence, Rhode Island, (2009).

\bibitem{Lieb1}
{\sc Lieb E., Loss M.} {\it Analysis.} Providence, Rhode Island, 1997.

\bibitem{Ostrovsky1}
{\sc  Ostrovsky E.I.} (1999). {\it Exponential estimations for random Fields and its
Applications, (in Russian).}  Moscow - Obninsk, OINPE.

\bibitem{Ostrovsky2}
{\sc Ostrovsky E.I.} (1980).{\it On the support of probabilistic measures in separable
Banach spaces.} Soviet Mathematic, Doklady, v.255, No 6 pp. 836 - 838, (in Russian).

\bibitem{Osekowski1}
{\sc A. Osekowski.} {\it Inequalities for dominated martingales.} Bernoulli 13 (2007), 54-79.

\bibitem{Osekowski2}
{\sc A. Osekowski.}  {\it Sharp martingale and semimartingale inequalities.} Monografie Matematyczne 72,
Birkh¨auser, 2012.


\bibitem{Ostrovsky3}
{\sc  Ostrovsky E. and Sirota L.} {\it  Moment and tail estimates for martingales and martingale transform,
with application to the martingale limit theorem in Banach spaces.}
arXiv:1206.4964v1 [math.PR] 21 Jun 2012

\bibitem{Ostrovsky603}
{\sc Ostrovsky E.} {\it Support of Borelian measures in separable Banach spaces. }
arXiv:0808.3248v1 [math.FA] 24 Aug 2008

\bibitem{Ostrovsky502}
{\sc  Ostrovsky E. and Sirota L.}  {\it  Schl\"omilch and Bell series for Bessel's functions, with
probabilistic applications.}
arXiv:0804.0089v1 [math.CV] 1 Apr 2008

\bibitem{Talagrand1}
 {\sc Talagrand M.} (1996). {\it Majorizing measure: The generic chaining.}
 Ann. Probab., {\bf 24} 1049 - 1103. MR1825156

\bibitem{Talagrand2}
 {\sc Talagrand M.} (2005). {\it The Generic Chaining. Upper and
     Lower Bounds of Stochastic Processes.} Springer, Berlin. MR2133757.

\vspace{16mm}

\bibitem{Araujo1}
{\sc Araujo  A., Gine E. }
{\it The central limit theorem for real and Banach valued random variables.}
Wiley, (1980), London, New York.

\bibitem{Billingsley1}
{\sc Billingsley P.} {\it Probability and measure.}
Wiley, 1979, London, New York.

\bibitem{Billingsley2}
{\sc Billingsley P.} {\it  Convergence of probability measures.}
Wiley, (1968), London, New York.

\bibitem{Dudley1}
{\sc Dudley R.M.} {\it Uniform Central Limit Theorem}. Cambridge University Press, (1999)

\bibitem{Grenander1}
{\sc Grenander U.} {\it Probabilities on algebraic structures.}
Wiley, 1963; London, New York.

\bibitem{Ledoux1}
 {\sc Ledoux M., Talagrand M.} (1991) {\it Probability in Banach Spaces.}
      Springer, Berlin, MR 1102015.

\bibitem{Ostrovsky1}
{\sc  Ostrovsky E.I.} (1999). {\it Exponential estimations for random Fields and its
applications (in Russian).}  Moscow - Obninsk, OINPE.

\vspace{16mm}

\bibitem{Fortet1}
{\sc Fortet R. and Mourier E.} {\it Les  fonctions  alratoires comme elements aleatoires dans
les espaces de Banach.} Studia Math., {\bf 15}, (1955), 62-79.

\bibitem{Garling1}
{\sc Garling D.J.H. }
{\it Functional Central Limit Theorems in Banach Spaces.}
 The Annals of Probability, Vol. 4, No. 4 (Aug., 1976), pp. 600-611

\bibitem{Gine1}
{\sc Gine E.} {\it On the Central Limit theorem for sample continuous processes.} Ann.
Probab. (1974), 2, 629-641.

\bibitem{Gine2}
{\sc Gine E., Zinn J.} {\it  Central Limit Theorem and Weak Laws of Large Numbers in certain Banach Spaces.  }
Z. Wahrscheinlichkeitstheory verw. Gebiete. {\bf 62}, (1983), 323  -  354.

\bibitem{Heinkel1}
{\sc Heinkel B.} Measures majorantes et le theoreme de la limite centrale dans
C(S). Z. Wahrscheinlichkeitstheory. verw. Geb., (1977). 38, 339-351.

\bibitem{Jain1}
{\sc Jain N.C. and Marcus M.B.} {\it Central limit theorem for $C(S)$ valued random
variables.} J. of Funct. Anal., (1975), 19, 216-231.

\bibitem{Kozachenko1}
 {\sc Kozachenko Yu. V., Ostrovsky E.I.} (1985). {\it The Banach Spaces of
      random Variables of subgaussian type.} Theory of Probab. and Math.
      Stat. (in Russian). Kiev, KSU, {\bf 32}, 43 - 57.

\bibitem{Ostrovsky301}
{\sc Ostrovsky E., L.Sirota L.}
{\it CLT for continuous random processes under approximations terms.}
arXiv:1304.0250v1 [math.PR] 31 Mar 2013

\bibitem{Pisier G.}
{\sc G. Pisier, J. Zinn.} {\it On the limit theorems for random variables with values in the spaces } $ L_p, \ 2  \le p  < \infty. $
 Z. Wahrscheinlichkeitstheorie verw. Gebiete 41, 289 = 304 (1978).

\bibitem{Prokhorov1}
{\sc Prokhorov Yu.V.} {\it Convergense of Random Processes and Limit Theorems
of Probability Theory.} Probab. Theory Appl., (1956), V. 1, 177-238.

\bibitem{Rackauskas1}
{\sc Rackauskas A, Suquet Ch.}
{\it Central limit theorems in H\"ölder topologies for Banach space valued random fields.}
Teor. Veroyatnost. i Primenen., 2004, Volume 49, Issue 1, Pages 109–125 (Mi tvp238)

\bibitem{Rhee1}
{\sc Rhee WanSoo and Michel Talagrand M. }
{\it Uniform bound in the central limit theorem for Banach space valued dependent random variables},
 Journal of Multivariate Analysis, 1986, vol. 20, issue 2, pages 303-320

\bibitem{Song1}
{\sc Song L.} {\it A counterexample in the Central Limit Theorem}.
Bulletin of the London Mathematical SocietyBulletin of the London Mathematical Society / Volume 31 /,
 Issue 02 / March 1999, pp 222-230

\bibitem{Sualb1}
{\sc Sualb Z.} {\it Central limit theorems for random processes with
sample paths in exponential Orlicz spaces.}
 Stochastic Processes and their Applications 66, (1997), l-20.

\bibitem{Zinn1}
{\sc Zinn J.  } {\it A Note on the Central Limit Theorem in Banach Spaces.}
 Ann. Probab. Volume 5, Number 2 (1977), 283-286.

\vspace{16mm}

\bibitem{Frolov1}
{\sc Frolov A.S., Tchentzov N.N. } {\it On the calculation by the Monte-Carlo
method definite integrals depending on the parameters. } Journal of Computetional
Mathematics and Mathematical Physics, (1962), V. 2, Issue 4, p. 714-718 (in
Russian).

\bibitem{Grigorjeva1}
{\sc Grigorjeva M.L., Ostrovsky E.I.} {\it Calculation of Integrals on discontinuous
Functions by means of depending trials method.} Journal of Computational
Mathematics and Mathematical Physics, (1996), V. 36, Issue 12, p. 28-39 (in
Russian).

\bibitem{Ostrovsky302}
{\sc Ostrovsky E., Sirota L.} {\it Monte-Carlo method for multiple parametric integrals
calculation and solving of linear integral Fredholm equations of a second
kind, with confidence regions in uniform norm.}
 arXiv:1101.5381v1 [math.FA] 27 Jan 2011

\bibitem{Ostrovsky303}
{\sc  Ostrovsky E., Rogover E.} {\it Non - asymptotic exponential bounds for
MLE deviation under minimal conditions via classical and generic chaining methods.}
arXiv:0903.4062v1 [math.PR] 24 Mar 2009

\vspace{4mm}

\end{thebibliography}
\end{document}